\documentclass[12pt,a4paper,twoside]{article}
\usepackage{amssymb,amsthm}
\usepackage[namelimits,sumlimits,fleqn]{amsmath}
\usepackage{setspace}
\usepackage[ hmargin={2.2cm,2.4cm}, top=36mm,
    a4paper,pdftex, textheight=230mm,
    headheight=20mm, headsep=10mm, bindingoffset=0.9cm]
    {geometry}
\usepackage[small, margin=20pt]{caption}
\usepackage{float}\restylefloat{figure}
\usepackage{url,hyperref}

\usepackage{fancyhdr}
\allowdisplaybreaks[4]

\title{New Proofs of Pl\"unnecke-type Estimates for Product Sets in Groups}
\author{Giorgis Petridis}
\date{}

\theoremstyle{plain}
\newtheorem{theorem}{Theorem}[section]
\newtheorem{lemma}[theorem]{Lemma}
\newtheorem{proposition}[theorem]{Proposition}
\newtheorem{corollary}[theorem]{Corollary}

\theoremstyle{definition}

\newtheorem*{acknowledgement}{Acknowledgement}

\theoremstyle{definition}

\newcommand{\rst}[1]{\ensuremath{{\mathbin\upharpoonright}%
\raise-.5ex\hbox{$#1$}}} 

\begin{document}

\onehalfspacing

\pagenumbering{arabic}

\setcounter{section}{0}

\bibliographystyle{plain}

\maketitle

\begin{abstract}
We present a new method to bound the cardinality of triple product
sets in groups and give three applications. A new and unexpectedly
short proof of the Pl\"unnecke-Ruzsa sumset inequalities for
Abelian groups. A new proof of a theorem of Tao on triple
products, which generalises these inequalities when no assumption
on commutativity is made. A further generalisation of the
Pl\"unnecke-Ruzsa inequalities in general groups.
\end{abstract}

\section[Introduction]{Introduction}
\label{Introduction}

Bounding the cardinality of sumsets is a central problem in
additive number theory and has many important applications. The
first bounds were obtained by Helmut Pl\"unnecke over forty years
ago in \cite{Plunnecke1969,Plunnecke1970}. Pl\"unnecke was
interested in the integers, but his graph-theoretic method works
equally well in any Abelian group. The upper bound he obtained is
still the standard:
\begin{theorem}[Pl\"unnecke]\label{Sumset Plunnecke}
Let $A$ and $B$ be finite sets in an Abelian group. Suppose
$$|A+B| \leq \alpha |A|.$$ Then there exists $X\subseteq A$ such that
$$|X+hB| \leq \alpha^h |X|.$$ In particular $|A+A|\leq \alpha |A|$
implies $|hA|\leq \alpha^h |A|$.
\end{theorem}
The bound is sharp in terms of $\alpha$ and $|X|$. Let for example
$A$ be a subgroup and $B$ a collection of generic points lying in
distinct cosets of $A$. Then $|A+B|=|A| |B|$ and so $\alpha=|B|$.
On the other hand $|A+hB| = \tbinom{|B|+h-1}{h} |X|$ for all
$X\subseteq A$.

Imre Ruzsa, who rediscovered Theorem \ref{Sumset Plunnecke}
\cite{Ruzsa1989,Ruzsa1991,Malouf1995}, extended Pl\"unnecke's
result to sum-and-difference sets. Combining Pl\"unnecke's theorem
with his `triangle inequality' (c.f Section \ref{Tao}) Ruzsa
showed:
\begin{theorem}[Ruzsa]\label{RuzsakB-lB}
Let $A$ and $B$ be finite sets in an Abelian group. Suppose that
$|A+B|\leq \alpha |A|$. Then $$|kB-lB|\leq \alpha^{k+l}|A|$$
provided that $k+l>1$.
\end{theorem}
The proofs given by Pl\"unnecke and Ruzsa were graph-theoretic and
to even provide a sketch one must start with a series of
definitions. Terence Tao obtained a purely combinatorial proof,
albeit with slightly worse bounds \cite{Granville2006,Tao-Vu2006}.
Under the assumption that $|A+A|\leq \alpha |A|$ he demonstrated
$$|kA-lA|\leq \alpha^{6(k+l)}|A|.$$ The bound may no longer be
sharp in $\alpha$, but the difference is not significant enough to
affect applications.

In recent years the study of product sets in not necessarily
Abelian groups has gained popularity. The lack of commutativity
imposes further restrictions, which make the general outlook quite
different. For example Theorem \ref{Sumset Plunnecke} no longer
holds. A well known counter example is $A=H\cup\{x\}$, where $H$
is a subgroup and $x$ is such that $|H x H|=|H|^2$. Then $|A
A|\leq 3 |A|$, while $|AAA|\geq (|A|-1)^2$. Ruzsa was probably the
first to realise \cite{Ruzsa2006} that an extra condition is
necessary should one aspire to extend Theorem \ref{Sumset
Plunnecke} to the non-Abelian case. In addition to $|AB|\leq
\alpha |A|$ one has to at least assume that
$$|A b B| \leq \alpha |A|$$ for all $b\in B$. Tao generalised the
above results to product sets for the important special case when
$A=B$ \cite{Tao2008}.
\begin{theorem}[Tao]\label{Tao Triple}
Let $B$ be a finite set in a group. Suppose that $|B b B| \leq
\alpha |B|$ for all $b\in B$ and also that $|B B| \leq \alpha
|B|$. Then $$|B^h| \leq \alpha^{ch}|B|$$ for some absolute
constant $c$.
\end{theorem}
The constant is large and is not calculated in \cite{Tao2008}.
Ruzsa asked in \cite{Ruzsa2010} for an explicit value of $c$ to be
calculated. It is well known \cite{Helfgott2008,Ruzsa2010,Tao2008} that it is
adequate to obtain a bound of the form $|BBB|\leq
\alpha^{c/2}|B|$.

Ruzsa proved a different extension to Theorem \ref{Sumset
Plunnecke} for $h=2$ by changing the order of multiplication and
focusing on $B A B$. By a clever application of Pl\"unnecke's
graph-theoretic method he showed in \cite{Ruzsa2009} the
following.
\begin{theorem}[Ruzsa]\label{Ruzsa Middle}
Let $A$, $B$ and $C$ be finite sets in a group. Suppose that $|A
B| \leq \alpha_1 |A|$ and that $|C A| \leq \alpha_2 |A|$. Then
there exists $X\subseteq A$ such that
$$|C X B| \leq \alpha_1 \alpha_2 |X|.$$
\end{theorem}
The example given above indicates that Theorem \ref{Ruzsa Middle}
does not lead to Pl\"unnecke-type bounds on the cardinality of
higher product sets $|AB^h|$.

A more thorough discussion of product set estimates in non-Abelian
groups can be found in \cite{Ruzsa2010,Tao2008}. In this note we
introduce a new and simpler method to obtain (variations of) the
above mentioned results.

In Section \ref{Ruzsa} we establish a variant of Ruzsa's theorem.
\begin{theorem}\label{Middle}
Let $A$ and $B$ be finite sets in a group. Suppose that $|A B|
\leq \alpha |A|$. Then there exists $X\subseteq A$ such that for
all finite sets $C$
$$|C X B| \leq \alpha |C X|.$$
\end{theorem}
Theorem \ref{Middle} is not as efficient as Theorem \ref{Ruzsa
Middle}. It only gives $|CXB|\leq \alpha_1 \alpha_2 |A|$. In
practice however one often bounds $|X|$ by $|A|$ and so the two
statements become equivalent. The advantage in using Theorem
\ref{Middle} is that the same $X$ works for all $C$, which helps
in some applications. Theorem \ref{Middle} furthermore has an
unexpectedly short proof despite the fact that we are multiplying
non-identical sets in a not necessarily Abelian group.

In Section \ref{Plunnecke} we deduce Theorems \ref{Sumset
Plunnecke} and \ref{RuzsakB-lB} from Theorem \ref{Middle}. The
deduction is swift and we thus present a short, elementary and
entirely self-contained proof, which results to the best known
bounds. It should also be noted that, by deducing Theorem
\ref{Sumset Plunnecke} from Theorem \ref{Middle}, we reverse the
usual order of doing things.

In Sections \ref{Tao} and \ref{NAP A comp to B} we study
non-Abelian analogues of Theorem \ref{Sumset Plunnecke}. In
section \ref{Tao} we materialise Ruzsa's suggestion and prove an
explicit form of Tao's theorem:
\begin{theorem}\label{NAP A=B}
Let $B$ be finite a set in a group. Suppose that $|BB|\leq
\alpha|B|$ and $|BbB|\leq \beta |B|$ for all $b\in B$. Then for
all $h>2$ $$ |B^h| \leq \alpha^{8h-17} \beta^{h-2} |B|.$$
\end{theorem}
$c$ in Theorem \ref{Tao Triple} can thus be taken to be nine. Our
approach is inspired by Tao's paper, but by using Theorem
\ref{Middle} we get a better dependence on $\alpha$, $\beta$.

It is well know that a similar approach works in a more general
setting. Under the further assumption that $A$ and $B$ have
comparable sizes we establish in Section \ref{NAP A comp to B} a
further generalisation of Theorem \ref{Sumset Plunnecke} to the
non-Abelian setting.
\begin{theorem}\label{NAP A<B}
Let $A$ and $B$ be finite sets in a group. Suppose that
\begin{itemize}
\item[(1)]$|AB|\leq \alpha |A|$. \item[(2)]$|AbB|\leq \beta |A|$
for all $b\in B$. \item[(3)] $|A|\leq \gamma|B|$.
\end{itemize}
Then there exists $S\subseteq A$ such that for all $h>1$ $$|SB^h|
\leq\alpha^{8h-9} \beta^{h-1} \gamma^{4h-5} |S|.$$
\end{theorem}
\begin{acknowledgement}
The author would like to thank Tim Gowers for his suggestion to
look at the non-Abelian setting and other recommendations that
improved the presentation of this note.
\end{acknowledgement}

\section{Growth of Triple Products}
\label{Ruzsa}

Our method is based on the choice of the subset $X$. We chose it
to be a subset of $A$ that grows minimally under multiplication by
$B$. The motivation for doing this comes from Pl\"unnecke's
original graph-theoretic method. A more illuminative explanation
of why this is a natural choice can be found in
\cite{Gowers-NAP2011}. More specifically in the first comment by
Tim Gowers.

Let us begin by explaining what we mean by minimal growth under
multiplication by $B$. For any $Z\subseteq A$ we define the
quantity $$ r(Z) = \frac{|Z B|}{|Z|}.$$ We let
$$K=\min_{Z\subseteq A} r(Z)$$ so that $|ZB|\geq K |Z|$ for all
$Z\subseteq A$. We choose $X$ to be such that $r(X)=K$.
With this in mind we prove a slightly more technical result.
\begin{proposition}\label{Stronger Middle} Let $X$ and $B$ be
finite sets in a group. Suppose that
$$K:=\frac{|X B|}{|X|} \leq \frac{|Z B|}{|Z|}$$ for all $Z\subseteq X$.
Then for all finite sets $C$ $$ |CXB|\leq K |CX| =
\frac{|CX|\,|XB|}{|X|}.$$
\end{proposition}

\begin{proof}
Let $C= \{c_1,\dots,c_r\}$. Using this (arbitrary) order of the
elements of $C$ we write $$ CX = \bigcup_{i=1}^r (c_i X_i)$$ where
$X_1=X$ and for $i>1$ $$X_i = \{x\in X : c_i x \notin
\{c_1,\dots,c_{i-1}\}X \}.$$ Observe that for all $j$:
\begin{eqnarray*}
\{c_1,\dots,c_j\} X = \bigcup_{i=1}^j (c_i X) = \bigcup_{i=1}^j
(c_i X_i).
\end{eqnarray*}
The sets $c_iX_i$ are disjoint and so for all $j$:
\begin{eqnarray}\label{X_i property}
\left|\{c_1,\dots,c_j\} X\right| = \sum_{i=1}^j
|c_iX_i|=\sum_{i=1}^j |X_i|.
\end{eqnarray}
We proceed by induction on $r$. For $r=1$ we have
$|c_1XB|=|XB|=K|X|=K|c_1X|$. For $r>1$ let us write $X_r^c = X
\backslash X_r$ for the complement of $X_r$ in $X$. By the
definition of $X_r$ we have that $c_r X_r^c \subseteq
\{c_1,\dots,c_{r-1}\} X$ and thus $c_r X_r^c B \subseteq
\{c_1,\dots,c_{r-1}\} X B$. Hence
\begin{eqnarray*}
C X B = \{c_1,\dots,c_r\}X B \subseteq (\{c_1,\dots,c_{r-1}\}X B)
\cup ((c_r X B)\backslash(c_r X_r^cB)).
\end{eqnarray*}
Note that $|(c_rXB)\backslash(c_rX_r^cB)|=
|(XB)\backslash(X_r^cB)|=|XB|-|X_r^cB|$ and so in particular
\begin{eqnarray}\label{inductive sum}
|C X B| \leq |\{c_1,\dots,c_{r-1}\}X B| + (|X B|-|X_r^cB|).
\end{eqnarray}
The first summand in the expression \eqref{inductive sum} above is
bounded by the inductive hypothesis on $r$ and \eqref{X_i
property}:
\begin{eqnarray*}
|\{c_1,\dots,c_{r-1}\} X B| \leq K |\{c_1,\dots,c_{r-1}\}X| = K
\sum_{i=1}^{r-1}|X_i|.
\end{eqnarray*}
The second bracketed term in \eqref{inductive sum} is at most
$K|X_r|$ as
\begin{eqnarray*}
|X B|- |X_r^cB|            & \leq & K|X|- K|X_r^c|\\
                           &   =  & K (|X|-|X_r^c|)\\
                           &   =  & K |X_r|
\end{eqnarray*}
the inequality following from the condition in the statement of
the proposition. Adding these upper bounds for the two terms in
\eqref{inductive sum} gives
\begin{eqnarray*}
|C X B| \leq K \sum_{i=1}^r |X_i|
\end{eqnarray*}
and the proposition follows by \eqref{X_i property}.
\end{proof}
Proposition \ref{Stronger Middle} is best possible even in the
Abelian case as we see by taking $C$, $X$ and $B$ to be groups
thought of as sets in the Cartesian product $C\times X \times B$.
It should also be noted that, as Ruzsa observed in
\cite{Ruzsa2009}, Proposition \ref{Stronger Middle} is a somewhat
commutative result with associativity playing a crucial role.
Theorem \ref{Middle} follows immediately:

\begin{proof}[Proof of Theorem \ref{Middle}]
We choose $X\subseteq A$ such that $$\frac{|XB|}{|X|}\leq
\frac{|ZB|}{|Z|}$$ for all $Z\subseteq A$ and apply Proposition
\ref{Stronger Middle} observing that $$ K= \frac{|XB|}{|X|}\leq
\frac{|AB|}{|A|} \leq\alpha.\qedhere$$
\end{proof}

\section{Pl\"unnecke-Ruzsa Inequalities}
\label{Plunnecke}

We now apply Proposition \ref{Stronger Middle} repeatedly to
deduce a slightly stronger version of Theorem \ref{Sumset
Plunnecke} where the subset $X$ is the same for all $h$.
Commutativity is crucial, but in a subtle way.
\begin{theorem}\label{Stronger Sumset Plunnecke} Let $A$ and
$B$ be finite sets in an Abelian group. Suppose that
$$|A+B| \leq \alpha |A|.$$ Then there exists $X\subseteq A$ such
that $$|X+hB| \leq \alpha^h |X|$$ holds for all $h$.
\end{theorem}

\begin{proof}
This is done by induction on $h$. Let $X\subseteq A$ be such that
$$\frac{|X+B|}{|X|}\leq \frac{|Z+B|}{|Z|}$$ for all $Z\subseteq
A$. For $h=1$ simply observe $|X+B|\leq |X|\,|A+B|/|A|\leq \alpha
|X|$. For $h>1$ we let $C=(h-1)B$. The condition in the statement
of Proposition \ref{Stronger Middle} is satisfied and thus
$$|X+hB|=|(h-1)B+X+B|\leq \alpha |X+(h-1)B| \leq \alpha^h
|X|.\qedhere$$
\end{proof}
There are circumstances where distinguishing between $|X+B|/|X|$
and $\alpha$ is worthwhile (e.g. in \cite{GPCaSu}). In most cases
however taking $K=\alpha$, like we implicitly did, is adequate.
One can strengthen Pl\"unnecke's graph-theoretic inequality along
the lines of Theorem \ref{Stronger Sumset Plunnecke}. Details and
applications can be found in \cite{GPCaSu}.

Using this stronger form of Theorem \ref{Sumset Plunnecke}
simplifies slightly the proof of Theorem \ref{RuzsakB-lB}.
\begin{proof}[Proof of Theorem \ref{RuzsakB-lB}]
We apply the triangle inequality of Ruzsa \cite{Ruzsa1978} that
will also be used in the next section. Let $X$, $Y$ and $Z$ be
finite sets in an Abelian group. Then
\begin{eqnarray}\label{ruzsa triangle commutative}
|X|\,|Y-Z|\leq |X+Y|\,|X+Z|.
\end{eqnarray}
Setting $X=X$, $Y=kB$ and $Z=lB$ gives $$|X|\,|kB-lB| \leq
|X+kB|\,|X+lB| \leq \alpha^{k+l} |X|^2 \leq \alpha^{k+l} |X|
|A|.\qedhere$$
\end{proof}
In the traditional deduction it is not enough to apply Theorem
\ref{Plunnecke}, but rather the graph-theoretic inequality from
which it follows.

\section{The Non-Abelian Setting: Tao's Theorem}
\label{Tao}

We now turn to the non-Abelian case and prove Theorem \ref{NAP
A=B}. The material in this section is similar to Tao's argument in
\cite{Tao2008}. A key difference is that we use Proposition
\ref{Stronger Middle}, which should be thought of as a non-Abelian
analogue of Theorem \ref{Sumset Plunnecke}. This simplifies the
argument and also results in a better dependence on $\alpha$,
$\beta$.

The proofs require two results of Ruzsa. Ruzsa's covering lemma
\cite{Ruzsa1999}:
\begin{lemma}[Ruzsa]\label{Ruzsa Covering}
Let $A$ and $B$ be finite sets in a group. Suppose that $|A B|
\leq K |A|$. Then there exists a set $S\subseteq B$ of size at
most $K$ such that $B\subseteq A^{-1} A S$.
\end{lemma}
and Ruzsa's triangle inequality \cite{Ruzsa1978}, a non-Abelian
generalisation of \eqref{ruzsa triangle commutative} \cite{Helfgott2008,Tao2008} :
\begin{lemma}[Ruzsa]\label{Ruzsa Triangle}
Let $X$, $Y$ and $Z$ be finite sets in a group. Then
$$|X|\,|Y Z|\leq |Y X^{-1}|\,|X Z|.$$ 
\end{lemma}
In both \cite{Ruzsa2010,Tao2008} it is shown how the above lemma
allows one to pass from an upper bound on triple products to upper
bounds on higher products. So in principle it is adequate to get a
Pl\"unnecke-type bound for $|BBB|$. In practice it is a little
more efficient to combine Lemma \ref{Ruzsa Triangle} with
Proposition \ref{Stronger Middle}. We will need the following
calculation in both steps of the argument.
\begin{corollary}\label{B-AA-B}
Let $A$ and $B$ be finite sets in a group. Suppose that $|BB|\leq
\alpha |B|$ and that $|BAB|\leq \alpha^2 |B|$. Then
$|BA^{-1}AB^{-1}|\leq \alpha^6 |B|$.
\end{corollary}
\begin{proof}
We apply Lemma \ref{Ruzsa Triangle} with $X=B$, $Y=BA^{-1}$ and
$Z=AB^{-1}$ and get
\begin{eqnarray}\label{to B-AB}
|BA^{-1}AB^{-1}|\leq \frac{|BA^{-1}B^{-1}|\,|BAB^{-1}|}{|B|} =
\frac{|BAB^{-1}|^2}{|B|}.
\end{eqnarray}
To bound $BAB^{-1}$ we once again apply Lemma \ref{Ruzsa
Triangle}. This time we set $X=B^{-1}$, $Y=BA$ and $Z=B^{-1}$
\begin{eqnarray}\label{BA-B}
|BAB^{-1}|\leq \frac{|BAB|\,|B^{-1}B^{-1}|}{|B^{-1}|} =
\frac{|BAB|\,|BB|}{|B|}\leq \alpha^3 |B|.
\end{eqnarray}
Putting \eqref{to B-AB} and \eqref{BA-B} together finishes the
proof.
\end{proof}
We now prove that sets which satisfy the conditions of Theorem
\ref{Triple} have small tripling.
\begin{theorem}\label{Triple}
Let $B$ be a finite set in a group. Suppose that $|BB|\leq \alpha
|B|$ and that $|B b B|\leq \beta |B|$ for all $b\in B$. Then
$$|BBB|\leq \alpha^7\beta\, |B|.$$
\end{theorem}
The strategy of the proof is to first choose $A\subseteq B$ that
grows minimally under multiplication by $B$. Then use Lemma
\ref{Ruzsa Covering} to cover the middle copy of $B$ by
$A^{-1}AX$. Finally apply Lemma \ref{Ruzsa Triangle} repeatedly to
break down the product in terms that are easily bound.
\begin{proof}[Proof of theorem \ref{Triple}]
We select $A\subseteq B$ such that $$ K:= \frac{|AB|}{|A|}\leq
\frac{|ZB|}{|Z|}$$ for all $Z\subseteq B$ and observe that $K\leq
\alpha$. Applying Proposition \ref{Stronger Middle} for $C=B$
gives
\begin{eqnarray}\label{BAB}
|BAB| \leq K |BA| \leq \alpha |BB| \leq \alpha^2 |B|.
\end{eqnarray}
Next we apply Lemma \ref{Ruzsa Covering} and get $T\subseteq B$ of
size at most $\alpha$ such that $B\subseteq A^{-1} A T$. In
particular we have
$$BBB \subseteq B A^{-1} A T B.$$ By setting $X=B$, $Y=B A^{-1} A$
and $Z=T B$ in Lemma \ref{Ruzsa Triangle} we get
\begin{eqnarray}\label{3 to 4}
|B|\,|BBB|\leq |B|\,|BA^{-1}ATB|\leq |BA^{-1}AB^{-1}|\,|BTB| .
\end{eqnarray}
To bound the second term in \eqref{3 to 4} we use the second
condition on the theorem together with the fact that $T\subseteq
B$.
\begin{eqnarray}\label{2nd term}
|B T B| = \left|\bigcup_{t\in T} (B t B) \right| \leq \sum_{t\in
T} |B t B| \leq |T| \beta |B| \leq \alpha\beta\,|B|.
\end{eqnarray}
The first term in \eqref{3 to 4} is at most $\alpha^6 |B|$ as we
see by combining \eqref{BAB} and Corollary \ref{B-AA-B}.
Substituting this in \eqref{3 to 4} and using \eqref{2nd term}
gives the desired bound on $|BBB|$.
\end{proof}
A slightly weaker bound for $|BBB|$ can be obtained using Theorem
\ref{Ruzsa Middle}. In this case we do not know whether $|AB|\leq
\alpha|A|$ and so have to content with $|AB|\leq |BAB|\leq
\alpha^2 |A|$ -- the first inequality following from the fact that
for any $b\in B$ we have $|AB|=|bAB|$. This leads to $|T|\leq
\alpha^2$ and we therefore only get $|BBB|\leq \alpha^9 \beta$.

Theorem \ref{NAP A=B} is completed by induction on $h$.
\begin{proof}[Proof of Theorem \ref{NAP A=B}] Assume that $h>2$
and select $A\subseteq B$ such that $$ K:= \frac{|AB|}{|A|}\leq
\frac{|ZB|}{|Z|}$$ for all $Z\subseteq B$ and observe that $K\leq
\alpha$. Applying Proposition \ref{Stronger Middle} for $C=B$
gives \eqref{BAB}. Lemma \ref{Ruzsa Covering} implies that there
exists $T\subseteq B$ of size at most $\alpha$ such that
$B\subseteq A^{-1} A T$. In particular
$$B^h \subseteq B A^{-1} A T B^{h-2}.$$
By setting $X=B$, $Y=BA^{-1}A$ and $Z=TB^{h-2}$ we get
\begin{eqnarray}\label{Bh to Bh-1}
|B^h|\leq |BA^{-1}ATB^{h-2}|\leq \frac{|BA^{-1}AB^{-1}|}{|B|}
|BTB^{h-2}|\leq \alpha^6 |BTB^{h-2}|.
\end{eqnarray}
The final inequality comes from \eqref{BAB} and Corollary
\ref{B-AA-B}. To bound $|BTB^{h-2}|$ we apply Lemma \ref{Ruzsa
Triangle} with $X=B^{-1}$, $Y=BT$ and $Z=B^{h-2}$.
\begin{eqnarray}\label{BTBh-2}
|BTB^{h-2}| \leq \frac{|BTB|}{|B^{-1}|} |B^{-1}B^{h-2}|\leq \alpha
\beta |B^{-1}B^{h-2}|.
\end{eqnarray}
The second inequality comes from \eqref{2nd term}. To bound
$|B^{-1}B^{h-2}|$ we once again resort to Lemma \ref{Ruzsa
Triangle} this time setting $X=B$, $Y=B{-1}$ and $Z=B^{h-2}$.
\begin{eqnarray}\label{-BBh-2}
|B^{-1}B^{h-2}| \leq \frac{|B^{-1}B^{-1}|}{|B|} |B^{h-1}| \leq
\frac{|BB|}{|B|} |B^{h-1}| \leq \alpha |B^{h-1}|.
\end{eqnarray} Putting \eqref{Bh to Bh-1}, \eqref{BTBh-2} and
\eqref{-BBh-2} gives $$ |B^h| \leq \alpha^8 \beta |B^{h-1}|$$ and
the desired bound on $|B^h|$ follows by the inductive hypothesis.
\end{proof}

As Ruzsa observed in \cite{Ruzsa2010} Theorem \ref{Triple} leads,
via a slightly different application of Lemma \ref{Ruzsa
Triangle}, to:
\begin{theorem}\label{NAP A=B Alternate Products}
Let $B$ be a finite set in a group and let
$\epsilon_1,\dots,\epsilon_h \in \{+1,-1\}$. Suppose that
$|BB|\leq \alpha|B|$ and $|BbB|\leq \beta |B|$ for all $b\in B$.
Then $$ |B B^{\epsilon_1}\cdots B^{\epsilon_h} B^{-1}| \leq
(\alpha^7 \beta)^{2h} |B|.$$
\end{theorem}

\section{The Non-Abelian Setting: Proof of Theorem \ref{NAP A<B}}
\label{NAP A comp to B}

Theorem \ref{NAP A<B} is proved similarly. We begin with the
result corresponding to Corollary \ref{B-AA-B}.
\begin{corollary}\label{S-SS-S}
Let $S$ and $B$ be finite sets in a group. Suppose that $|CSB|\leq
\alpha |CS|$ for all finite sets $C$ in the group. Then
$$|SS^{-1}SS^{-1}|\leq \alpha^6 \left(\frac{|S|}{|B|}\right)^3 |S|.$$
\end{corollary}
\begin{proof}
We begin by taking $X=B^{-1}$, $Y=SS^{-1}S$ and $Z=S^{-1}$ in
Lemma \ref{Ruzsa Triangle}
\begin{eqnarray*}
|SS^{-1}SS^{-1}| \leq \frac{|SS^{-1}SB|\,|B^{-1}S^{-1}|}{|B^{-1}|}
\leq \frac{\alpha |SS^{-1}S|\,|SB|}{|B|} \leq \alpha
\frac{|SS^{-1}S|\alpha|S|}{|B|}=\alpha^2 \frac{|SS^{-1}S|
|S|}{|B|}.
\end{eqnarray*}
The inequalities follow from the condition in the statement of the
corollary. For example $|SB|=|sSB|$ for any $s\in S$ and hence
$|SB|\leq \alpha |sS|=\alpha |S|$. To bound $|SS^{-1}S|$ we apply
Lemma \ref{Ruzsa Triangle} with $X=B^{-1}$, $Y=S$ and $Z=S^{-1}S$.
\begin{eqnarray*}
|SS^{-1}S| \leq \frac{|SB|\,|B^{-1}S^{-1}S|}{|B^{-1}|} \leq
\frac{\alpha|S|\,|S^{-1}SB|}{|B|} \leq \alpha
\frac{|S|\alpha|S^{-1}S|}{|B|} = \alpha^2
\frac{|S||SS^{-1}|}{|B|}.
\end{eqnarray*}
We are finally left with bounding $SS^{-1}$. We once again apply
Lemma \ref{Ruzsa Triangle} with $X=B^{-1}$, $Y=S$, and $Z=S^{-1}$.
\begin{eqnarray*}
|SS^{-1}| \leq \frac{|SB|\,|B^{-1}S^{-1}|}{|B^{-1}|} =
\frac{|SB|^2}{|B|} \leq \alpha^2 \frac{|S|^2}{|B|}.
\end{eqnarray*}
Putting everything together gives the desired bound.
\end{proof}
We next show that for any $S\subseteq A$, which grows minimally
under multiplication by $B$, $SBB$ satisfies a Pl\"unnecke-type
bound.
\begin{proposition}\label{SBB}
Let $A$ and $B$ be sets in a finite group. Suppose that
\begin{itemize} \item[(1)]$|AB|\leq \alpha |A|$.
\item[(2)]$|AbB|\leq \beta |A|$ for all $b\in B$. \item[(3)]
$|A|\leq \gamma|B|$.
\end{itemize}
Let $S\subseteq A$ be such that
$$\frac{|SB|}{|S|}\leq\frac{|ZB|}{|Z|}$$ for all $Z\subseteq A$.
Then $$ |SBB| \leq \alpha^7 \beta \gamma^3 |S|.$$
\end{proposition}
\begin{proof}
We begin by observing that Proposition \ref{Stronger Middle} can
be applied to $S$ and $B$ and so
\begin{eqnarray}\label{CSB}
|CSB|\leq \alpha |CS|.
\end{eqnarray}
In particular $|SB|\leq \alpha|S|$ and so by Lemma \ref{Ruzsa
Covering} there exists $T\subseteq B$ of size at most $\alpha$
such that $B\subseteq S^{-1}ST$. Thus $|SBB|\leq |SS^{-1}STB|$.
Applying Lemma \ref{Ruzsa Triangle} with $X=S$, $Y=SS^{-1}S$ and
$Z=TB$ gives
\begin{eqnarray}\label{SBB'}
|SBB|\leq |SS^{-1}STB|\leq\frac{|SS^{-1}SS^{-1}|}{|S|} |STB|.
\end{eqnarray}
By \eqref{CSB} and Corollary \ref{S-SS-S} we know that the first
term in \eqref{SBB'} is at most $\alpha^6 \gamma^2 |S|/|B|$. The
second term in \eqref{SBB'} is bounded using the second condition
and the fact that $S\subseteq A$ and $T\subseteq B$.
\begin{eqnarray}\label{second term}
|STB| = \left| \bigcup_{t\in T} (StB)\right| \leq \sum_{t\in T}
|StB| \leq \sum_{t\in T} |AtB| \leq |T| \beta |A| \leq
\alpha\beta\, |A|.
\end{eqnarray}
Substituting in \eqref{SBB'} gives $$|SBB|\leq \alpha^6 \gamma^2
\frac{|S|}{|B|}\,\alpha \beta |A|\leq \alpha^7 \beta \gamma^3 |S|.
\qedhere$$
\end{proof}
Theorem \ref{NAP A<B} follows by induction on $h$.
\begin{proof}[Proof of Theorem \ref{NAP A<B}]
We work with $S\subseteq A$ defined in the statement of
Proposition \ref{SBB}. Proposition \ref{Stronger Middle} implies
that \eqref{CSB} holds. Next we apply Lemma \ref{Ruzsa Covering}
and get $T\subseteq B$ of size at most $\alpha$ such that
$B\subseteq A^{-1} A T$. In particular
$$|SB^h| \subseteq |S S^{-1} S T B^{h-1}|.$$ By setting $X=S$, $Y=S S^{-1}
S$ and $Z=TB^{h-1}$ in Lemma \ref{Ruzsa Triangle} we get
\begin{eqnarray}\label{SBh}
|SB^h|\leq \frac{|SS^{-1}SS^{-1}|}{|S|} |STB^{h-1}|.
\end{eqnarray}
The second term is bounded by applying Lemma \ref{Ruzsa Triangle}
with $X=B^{-1}$, $Y=ST$ and $Z=B^{h-1}$.
\begin{eqnarray}\label{STBh-1}
|STB^{h-1}|\leq \frac{|STB| |B^{-1}B^{h-1}|}{|B^{-1}|}\leq \alpha
\beta \gamma |B^{-1}B^{h-1}|.
\end{eqnarray}
The second inequality follows from \eqref{second term}. We are
thus left to bound $|B^{-1}B^{h-1}|$. Setting $X=S$, $Y=B^{-1}$
and $Z=B^{h-1}$ in Lemma \ref{Ruzsa Triangle} gives
\begin{eqnarray}\label{-BBh-1}
|B^{-1}B^{h-1}|\leq \frac{|B^{-1}S^{-1}| |SB^{h-1}|}{|S|}=
\frac{|SB|}{|S|} |SB^{h-1}| \leq \alpha |SB^{h-1}|.
\end{eqnarray}
Substituting \eqref{STBh-1} and \eqref{-BBh-1} in \eqref{SBh} and
applying Corollary \ref{S-SS-S} we get
$$|SB^h| \leq \alpha^6 \gamma^3\,\alpha\beta\gamma\,\alpha |SB^{h-1}|=
\alpha^8 \beta \gamma^4 |SB^{h-1}|$$ and the theorem follows from
the inductive hypothesis.
\end{proof}

\bibliography{arxiv}

$\hspace{12pt}$\textsc{Department of Pure Mathematics and
Mathematical Statistics}\\ \textsc{Wilberforce Road, Cambridge CB3
0WB, England}

$\hspace{12pt}$ \textit{Email address}: giorgis@cantab.net

\end{document}